\documentclass[a4paper,12pt]{article}
\usepackage{amssymb}
\usepackage{amsmath}
\usepackage{graphicx}
\usepackage{epstopdf}

\newcommand{\eneqa}{\end{eqnarray}}
\newcommand{\begeqaet}{\begin{eqnarray*}}
\newcommand{\eneqaet}{\end{eqnarray*}}
\newcommand{\be}{\begin{equation}}
\newcommand{\ee}{\end{equation}}

\newcommand{\rn}{\rbig^n}

\newcommand{\rbig}{{\mathbb{R}}}

\newcommand{\opr}{\Omega^\prime}

\newcommand{\io}{\int_\Omega}
\newcommand{\iok}{I_{\mathrm{OK}}}

\newcommand{\ed}{\end{document}}

\def\beq{\begin{equation}}
\def\eeq{\end{equation}}

\def\Box{\hfill\framebox(0.25,0.25){}}

\newtheorem{thm}{Theorem}[section]

\newtheorem{defi}{Definition}[section]
\newtheorem{prop}{Proposition}[section]

\newtheorem{cor}{Corollary}[section]
\newcommand{\begeqa}{\begin{eqnarray}}

\begin{document}

\begin{center}{\bf\Large Global integrability and weak Harnack  estimates for elliptic PDE in divergence form}
\end{center}\smallskip

 \begin{center}
Boyan SIRAKOV\footnote{e-mail : bsirakov@mat.puc-rio.fr}\\
PUC-Rio, Departamento de Matematica,\\ Gavea, Rio de Janeiro - CEP 22451-900, BRAZIL
\end{center}\bigskip

{\small \noindent{\bf Abstract}. We show that two classically known properties of positive supersolutions of uniformly elliptic PDEs, the boundary point principle (Hopf lemma) and global integrability, can be quantified with respect to each other. We obtain an extension up to the boundary of the De Giorgi-Moser weak Harnack inequality, optimal with respect to the norms involved, for equations in divergence form.  }
%

\section{Introduction and Main Results}

This paper is devoted to global estimates for nonnegative supersolutions of divergence-form uniformly elliptic PDE in a given domain. We study bounds in terms of the distance to the boundary, as well as integrability and $L^p$-estimates up to the boundary, of supersolutions and their gradients.

A fundamental property of superharmonic functions is that positivity entails a quantitative version of itself.  Specifically, if $u>0$ is superharmonic in a bounded $C^{1,1}$-domain $\Omega\subset \rn$, then
\begin{equation}\label{hopf1}
u\ge c_0\, d\quad\mbox{ in }\;\Omega,\qquad \mbox{ where }\;d(x):=\mathrm{dist}(x,\partial\Omega),\;\;c_0= c_0(u,\Omega)>0.
\end{equation}
In other words, if $u$ attains a minimum at a boundary point, then the normal derivative of $u$ at this point does not vanish. This is the famous Zaremba-Hopf-Oleinik lemma, also called boundary point lemma or boundary point principle. Because of the importance of this principle, a lot of work has been dedicated to understanding its ramifications and getting optimal conditions for its validity, in terms of the regularity or the geometry of the domain, or of the nature of the coefficients of the elliptic operator. We refer to  \cite{AM}, \cite{AN1}, \cite{AN2}, \cite{CLN}, \cite{K1}, \cite{N1}, \cite{PS},  \cite{S1},  \cite{S2} where such conditions, as well as a lot more references  and  history of this ``bedrock" (to quote page 1 of \cite{PS}) result in the theory of elliptic PDE can be found.

Another striking property of superharmonic functions, to which a lot of attention has been given, is that positivity implies global integrability. A classical result by Armitage \cite{Ar1}, \cite{Ar2}, states that if $u>0$ is superharmonic in $\Omega$ then $u\in L^p(\Omega)$ for each $p<n/(n-1)$, and that bound is sharp. Extensions of Armitage's result to superharmonic functions in more general domains can be found in \cite{Ai}, \cite{MS}, \cite{SU}.  \medskip

The essence of the results below is that  {\it global integrability and the boundary point principle  quantify each other}.
Furthermore, we study and quantify how the loss of superharmonicity influences these properties - the integrability is preserved, with the boundary point estimate being corrected with  a $L^q$-norm of the ``loss", for $q>n$.

Our main result, Theorem \ref{bwhistrong} below,  provides a sharp global integrability estimate for the quantity $u/d$, and  can also be seen as an optimal global weak Harnack inequality for this quantity. A simpler variant of Theorem \ref{bwhistrong} is a {\it global extension of the classical De Giorgi-Moser weak Harnack inequality}, stated in Theorem \ref{bwhiweak}.

 Another  consequence of Theorem \ref{bwhistrong} is a novel and surprising global integrability result for the { gradient} of supersolutions, also quantified in terms of the boundary point property.\medskip

 We  consider general linear operators in divergence form, and weak solutions of inequalities in the form
\begin{equation}\label{first}
-Lu=-\mathrm{div}\left( A(x)\nabla u \right) + b(x) |\nabla u|\ge -f(x), \qquad u\ge 0,
\end{equation}
where $A$ is a symmetric matrix, for some $\lambda>0$ and $q>n$
\begin{equation}\label{hypocoef}
A\ge \lambda I,\qquad A\in W^{1,q}(\Omega),\qquad b,f\in L^q_+(\Omega),
\end{equation}
 $\Omega\subset \rn$, $n\ge2$, is a bounded $C^{1,1}$-domain. We set $\Lambda=\|A\|_{W^{1,q}(\Omega)}$. In the sequel, all constants denoted by $C$ will be allowed to depend  on $n$, $\lambda$, $\Lambda$, $q$, $\|b\|_{L^q(\Omega)}$,  the diameter of~$\Omega$, the $C^{1,1}$-norm of $\partial\Omega$, as well as on the positive exponents $p,s,$ in each of the theorems below.

\begin{defi}\label{def1} We say that $u:\Omega\to\mathbb{R}$ is a solution of \eqref{first} if for each $l\in \mathbb{N} $ the function $u_l:=\min\{u,l\}$ belongs to $H^1_{\mathrm{loc}}(\Omega)$, and
\begin{equation}\label{defsupersol}
\int_\Omega A\nabla u_l.\nabla \varphi +\int_\Omega\varphi b|\nabla u_l| \ge -\int_\Omega f\varphi,\qquad \mbox{for\ each}\; \varphi\in C^\infty_0(\Omega), \varphi\ge0.
\end{equation}
\end{defi}

In the literature there are at least four frequently used notions of supersolutions which we briefly recall: in the weak Sobolev sense ($u\in H^1_{\mathrm{loc}}(\Omega)$ and \eqref{defsupersol} holds for $u$ instead of $u_l$), in the $L^q$-viscosity sense ($u$ is continuous and  \eqref{first} holds in the essliminf sense for $W^{2,q}$-functions at points where they touch $u$ from below), in the $C$-viscosity sense (if  $A$, $b$, $f$ are continuous, $u$ is lower semi-continuous and  \eqref{first} holds for smooth functions at points where they touch $u$ from below), and in the potential theory sense ($u$ is continuous and is above the solution of the Dirichlet problem in any ball, with $u$ as boundary value). All these four definitions are included in Definition \ref{def1}, see the Appendix. We note that it is important to consider not just weak Sobolev supersolutions, in order to accommodate supersolutions which are not in the energy space $H^1_{\mathrm{loc}}$, such as the fundamental solution with pole inside the domain.

Notice that in Definition \ref{def1} we only ask that $u_l$ be in $H^1_{\mathrm{loc}}(\Omega)$, there is no a priori assumption on integrability or behaviour close to the boundary, let alone any boundary condition. The integrability is a consequence of the supersolution property only, as the following shows.\medskip

Here is our main result.

\begin{thm}\label{bwhistrong} Let $u$ be a nonnegative solution of \eqref{first}. Then
 \begin{equation}\label{ineq1}
\left( \int_{\Omega} \left(\frac{u}{d}\right)^s\right)^{1/s}\le C\left(\inf_{\Omega} \frac{u}{d} +\|f\|_{L^q(\Omega)}\right),
 \end{equation}
for each
$$
s<1.
$$
\end{thm}

In \cite[Theorem 1.2]{Sir1}  we showed that \eqref{ineq1} holds for some small $s>0$, for viscosity supersolutions of more general fully nonlinear equations. As a consequence, in \cite[Theorem 1.4]{Sir1} we obtained a global Harnack inequality for nonnegative solutions of { inhomogenous} equations. Here we manage to upgrade the estimate from \cite{Sir1} to the optimal range $s<1$, thanks to the additional variational structure we have.

 To our knowledge, even the finiteness of the integral in the left-hand side of \eqref{ineq1} is proved here for the first time, for values of $s$ not close to zero.

As a simple consequence of the proof of Theorem \ref{bwhistrong} we get the following up-to-the-boundary extension of the classical weak Harnack inequality.

\begin{thm}\label{bwhiweak} Let $u$ be a nonnegative solution of \eqref{first}. Then
 \begin{equation}\label{ineq2}
\left( \int_{\Omega} u^p\right)^{1/p}\le C\left(\inf_{\Omega} \frac{u}{d} +\|f\|_{L^q(\Omega)}\right),
\end{equation}
for each
$$
p<\frac{n}{n-1}.
$$
\end{thm}

The interior version of Theorem \ref{bwhiweak}, when $\Omega$ in the integral and the infimum in \eqref{ineq2} is replaced by a compactly included subdomain, is the famous and fundamental De Giorgi-Moser weak Harnack inequality, which is known to hold for  $p<n/(n-2)$ (see for instance \cite[Theorem 8.18]{GT}). It is worth noticing that the interior result does not require any regularity of the leading coefficients of the operator.
For  results  with a small exponent $p$ and equations with homogeneous operators like the $m$-Laplacian, we refer to \cite{L1}, \cite{Sa} (finiteness of $\|u\|_{L^p(\Omega)}$), and \cite{BM}, \cite{CLN} (bounds for  $\|u\|_{L^p(\Omega^\prime)}$, $\Omega^\prime\subset\subset\Omega$).\medskip

  In combination with known inequalities for supersolutions, Theorem~\ref{bwhistrong} implies optimal gradient integrability  and a gradient bound for nonnegative supersolutions.

\begin{thm}\label{bwhigrad}  Let $u$ be a nonnegative solution of \eqref{first}. Then
 \begin{equation}\label{ineq3}
\left( \int_{\Omega} |\nabla u|^s\right)^{1/s}\le C\left(\inf_{\Omega} \frac{u}{d} +\|f\|_{L^q(\Omega)}\right),
\end{equation}
for each
$$
s<1\,.
$$
\end{thm}

 That the integral in \eqref{ineq3} is finite for bounded solutions (as opposed to supersolutions) and $s<1$, and for supersolutions and some $s$ close to zero was proved in \cite{KK} for $m$-homogeneous equations ($f=0$, $m>1$). No upper bound for the integral in \eqref{ineq3} was previously known. See also the remark at the end of Section \ref{sectproofs}, below.\medskip

We stress that \eqref{ineq1}-\eqref{ineq3} are valid for supersolutions, there is no need for $u$ to satisfy an equation or even a two-sided inequality. When solutions are considered, there is an extensive literature on boundary Harnack inequalities in recent years, including equations with a right-hand side -- see \cite{CCCC}, \cite{B}, \cite{S1}, \cite[Theorem 2.4]{SS}, \cite[Theorem 1.4]{Sir1}, \cite{AS}. \medskip

Further, it is  important that reading \eqref{ineq1}-\eqref{ineq3} ``from right to left"
is a quantification of the boundary point lemma (i.e. of $c_0(u,\Omega)$ in \eqref{hopf1}), as well as an extension of this lemma to inhomogeneous inequations.
\begin{cor}\label{maincor} Let $u$ be a nonnegative solution of \eqref{first}. Then
$$
u(x)\ge \Big( {C}^{-1}\max\left\{  \|u\|_{L^p(\Omega)}, \|\nabla u\|_{L^s(\Omega)}, \left\|\frac{u}{d}\right\|_{L^s(\Omega)}\right\} -\|f\|_{L^q(\Omega)}\Big)\:d(x)
$$
for every $x\in\Omega$, and
$s<1$, $p<\frac{n}{n-1}$.
\end{cor}

Corollary \ref{maincor} is new even for $f=0$. The importance of this type of quantification of the boundary point principle was first recognised in \cite[Lemma 3.2]{BC} and \cite[Lemma 1.6]{CLN}, and has received a lot of attention in the recent years, with  the already quoted extensions in \cite{Sir1}, \cite{BM},  and applications in \cite{BFM}, \cite{NS},  \cite{Sir2}, \cite{STZ}. In particular, in \cite[Theorem 2]{Sir2} we prove a uniform $L^\infty$-bound for positive solutions of nonlinear two-sided inequalities, achieving optimal growth of the nonlinearities  thanks to the optimality of the ranges for $p,s$ in Theorems \ref{bwhistrong}-\ref{bwhiweak} above.\medskip

The ranges for $p,s$ in the above theorems are indeed optimal -- taking $\Omega$ to have a flat part of its boundary in $\{x_n=0\}$ containing the origin, and the harmonic function $u=x_n/|x|^n$, we easily see that the theorems fail for  $s=1$, resp. $p=n/(n-1)$. Thus, our results show that general supersolutions (with a priori arbitrarily bad behaviour in $\overline{\Omega}$) are as integrable as the Poisson kernel, i.e. the fundamental solution in a half-space with an one-point singularity on the boundary. This property may not look natural at a first sight, since supersolutions are supposed to be larger than solutions. That an one-sided elliptic inequality should imply gradient control is even less intuitive and more surprising. \medskip

Some comments on the assumptions on the elliptic operator are in order.
 The Sobolev regularity and integrability assumptions we made in \eqref{hypocoef} are not far from optimal -- it is known that even the Hopf lemma may fail for some $A\in C(\overline{\Omega})\cap W^{1,n}(\Omega)$ and $b=0$, or for $A=I$ and some $b\in L^n(\Omega)$ (see for instance \cite{N1}, \cite{S1}, \cite{AN1}).
 Given that the boundary point principle is valid for $A\in C^\alpha(\Omega)$, $\alpha>0$ (or even for $A$ Dini continuous in $\Omega$, see \cite{AN2}), it is tempting to conjecture that our results are also valid under this hypothesis. It is however worth noticing that the quantitative nature of the above theorems represents, in a certain sense, "propagation of smallness" of the quantity $u/d$ from some points to the whole of the domain; whereas in the recent years the theory on propagation of smallness has developed considerably: we refer to \cite{LM} for one such direction. Those profound results are valid for solutions of div$(A(x)Du)=0$, for a symmetric Lipschitz $A$, but are known to fail for $A\in C^\alpha$, $\alpha<1$. We use the symmetry of $A$ and $A\in W^{1,q}$ in the proof -- to write \eqref{thisone} and the chains of inequalities that follow with the help of the Divergence Theorem. We note notwithstanding that the symmetry hypothesis on $A$ is moot at least for strong (i.e. in $W_{\mathrm{loc}}^{2,q}$) supersolutions, since then we can write the inequality in non-divergence form (up to changing $b$) and replace $A$ by $(A+A^T)/2$.

\medskip

Broadly speaking, the main idea of the proof of Theorem~\ref{bwhistrong} is to use a Moser-type iteration in order to upgrade the result from \cite[Theorem 1.2]{Sir1} to the optimal range $s<1$. We rely on recent embedding results and estimates for weighted Sobolev spaces from \cite{FMT}.

It turns out that implementing an iteration procedure is considerably more delicate at the boundary than in the interior of the domain. The test function we use contains a product of different powers of $u$ and the distance function, and these powers are varied independently at some steps, and together at others. Surprisingly, it is indispensable to track carefully the dependence in these powers of the constants in front of the integrals in order to realize even one step of the iteration (no such necessity appears in the proof of the interior estimate). At several moments this dependence suffices just barely to absorb bad terms into good terms, see for instance \eqref{fine1}-\eqref{fine2}; in this sense the estimate \eqref{ineq1} feels very ``exact". Also a sequence of cut-offs that get close to the boundary is necessary, together with a careful evaluation of their contribution. At the end of the proof we obtain a somewhat unusual recursively defined sequence of Lebesgue exponents, which converges to one.

\section{Proofs}\label{sectproofs}

It is sufficient to prove inequalities \eqref{ineq1}-\eqref{ineq3} with $u$ replaced by $u_l$ (see Definition \ref{def1}), and a constant $C$ independent of $l$. Indeed,  then the monotone convergence theorem implies that $u$ satisfies the same inequalities. In particular, $u$ can be assumed bounded, i.e. a usual weak Sobolev supersolution, provided $C$ is shown to be independent of $u$. Further, we observe that $u$ is lower semi-continuous (see the appendix), in particular $u$ attains its minimum on compacts. We also recall that the minimum of supersolutions is a supersolution.

Thanks to the boundary weak Harnack inequality in \cite{Sir1} we know that \eqref{ineq1}-\eqref{ineq2} are true if $p,s<\varepsilon_0$, where $\varepsilon_0$ is a small positive constant which depends on the right quantities. We note that \cite[Theorem 1.2]{Sir1} was stated for $L^q$-viscosity supersolutions but it also applies to bounded  supersolutions  as in Definition~\ref{def1} -- see the Appendix.

So our goal will be to improve \cite[Theorem 1.2]{Sir1} to every $\varepsilon_0<1$.

\subsection{Proof of Theorem \ref{bwhistrong}}

For readers' convenience  we start by giving a local and scaled version of Theorem \ref{bwhistrong}. We denote with $B_R^+=\{x\in\rn\::\:|x|<R,\;x_n>0\}$ a half-ball whose boundary's flat portion is included in $\{x_n=0\}$.

\begin{thm}\label{bwhistrongloc} Assume that $u_l\in H^1_{\mathrm{loc}}(B_{2R}^+)$ is a bounded weak Sobolev solution of \eqref{first} in $B_{2R}^+$, for each $l\in\mathbb{N}$. Then there exists $C>0$ depending on $n$, $\lambda$, $\Lambda$, $q$, $s$, and $R^{1-n/q}\|b\|_{L^q(B_{2R}^+)}$, such that
 \begin{equation}\label{ineq1loc}
R^{-n/s}\left( \int_{B_{3R/2}^+} \left(\frac{u}{x_n}\right)^s\right)^{1/s}\le C\left(\inf_{B_R^+} \frac{u}{x_n} +R^{1-n/q}\|f^-\|_{L^q(B_{2R}^+)}\right),
 \end{equation}
for each
$$
s<1.
$$
\end{thm}

By a standard argument involving local straightening and covering of $\partial \Omega$, Theorem \ref{bwhistrong}  is a  consequence of Theorem \ref{bwhistrongloc}. Note we do not strictly need to use Theorem \ref{bwhistrongloc} in order to prove  Theorem \ref{bwhistrong} (which is easily seen to be equivalent to Proposition \ref{propgen} below), however we include Theorem \ref{bwhistrongloc} since we believe it is useful to display a local and rescaled with respect to the size of a domain version of the main theorem. The main results  in \cite{Sir1} were also stated in half-balls.

The following proposition is  Theorem \ref{bwhistrong} for a bounded convex $C^2$-domain and $u\in H^1$ up to the boundary.
\begin{prop}\label{propgen}  Let  $\Omega$ be a bounded convex $C^2$-domain, and $u\in H^1(\Omega)$ be a bounded weak Sobolev solution of \eqref{first} in $\Omega$.  Then
 \begin{equation}\label{ineq1conv}
\left( \int_{\Omega} \left(\frac{u}{d}\right)^s\right)^{1/s}\le C\left(\inf_{\Omega} \frac{u}{d} +\|f\|_{L^q(\Omega)}\right),
 \end{equation}
for each $s<1$.
\end{prop}\medskip

\noindent {\it Proof of Theorem \ref{bwhistrongloc} assuming Proposition \ref{propgen}}. By scaling ($x\to x/R$), it is enough to prove \eqref{ineq1loc} for $R=1$. Fix a smooth convex domain $\Omega$, such that $  B_{4/3}^+\subset \Omega\subset B_2^+$, with the $C^2$-norm of $\partial \Omega$ being a universal constant, and take a monotone sequence of smooth convex domains $\omega_m$ which converges to $\Omega$ in $C^2$ and $\omega_m\subset \Omega\cap\{x_n>1/m\}$. Since $u\in H^1(\omega_m)$ we can apply \eqref{ineq1conv} with $\Omega$ replaced by $\omega_m$, and then pass to the limit $m\to\infty$ in the resulting inequalities,  with the help of the monotone convergence theorem.\hfill $\Box$
\medskip

The rest of this section  will be devoted to the proof of Proposition \ref{propgen}. \medskip

We are going to use the following weighted Sobolev inequalities, which follow from a result due to Filippas, Maz'ya and Tertikas, \cite{FMT}.
\begin{thm}[\cite{FMT}]\label{wsob} Let $\Omega$ be a bounded convex $C^2$-domain of $\mathbb{R}^n$, $n>2$. If $\phi\in H^1_0(\Omega)$ then for all
$$a\in\left(\frac{1}{2}\,,1\right),\qquad t\in \left(2, \frac{2n}{n-2}\right] $$
we have
$$
\|d^b\phi\|_{L^t(\Omega)}\le C\| d^{a}\nabla \phi\|_{L^2(\Omega)}  + C\|\phi\|_{L^2(\Omega)}
$$
where $C=C(n,t,\Omega)$, and we have set
$$
b= a-1 + \frac{t-2}{2t}\,n.
$$

More generally, if $n>\alpha>1$,
$$
\|d^b\phi\|_{L^t(\Omega)}\le C\| d^{a}\nabla \phi\|_{L^\alpha(\Omega)}  + C\|\phi\|_{L^\alpha(\Omega)},\qquad b= a-1 + \frac{t-\alpha}{\alpha t}\,n.
$$
\end{thm}

\noindent {\it Proof.} Theorem \ref{wsob} is a particular case of  \cite[Theorem 4.5]{FMT} and inequality (4.40) in \cite{FMT}\footnote{note we apply these results with $k=1$ in \cite{FMT}; and that there is a misprint in (4.40) in \cite{FMT}, the $L^q$-norm in the right-hand side lacks the power $p$ there.}.

Alternatively, Theorem \ref{wsob} is a consequence of  \cite[Proposition 2]{So}, combined with the H\"older inequality and the interpolation Lemma 3 in that paper (to accomodate the reader we note that $a$ in \cite{So} is our $2a$, while $b$ in \cite[Lemma 3]{So} is our $2-2a$).\hfill $\Box$ \bigskip

\noindent {\it Proof of  Proposition \ref{propgen}}.
We know there is a uniform neighborhood of size $\delta>0$ of the boundary $\partial\Omega$ in which the distance function to the boundary is $C^2$-smooth.  By scaling we can assume that $\delta=2$ (translate so that $0\in\Omega$, and dilate $x\to R_0x$, for some $R_0$ which depends on $\delta$, diam$(\Omega)$, and $\min_{x\in\partial\Omega}\cos(x,-\nu(x))>0$, where $\nu(x)$ is the interior normal to the boundary of~$\Omega$ at $x$). Set $\opr=\Omega_1^\prime=\{x\in\Omega\::\: \mathrm{dist}(x,\partial \Omega)<1\}$, $\Omega^{\prime\prime}=\{x\in\Omega\::\: \mathrm{dist}(x,\partial \Omega)<2\}$, and   $\Omega_m^\prime=\{x\in\Omega\::\: \mathrm{dist}(x,\partial \Omega)<1/m\}$, $m\in\mathbb{N}$, $m\ge1$.

We fix a $C^2$-smooth unitary vector field $\nu(x)$ in $\overline{\Omega}$, such that for each $x\in\opr$ with $d(x)=1/m$,  $\nu(x)$ is the interior normal to the boundary of $\Omega_m^\prime$ at $x$, and for $x\in\partial\Omega$, $\nu(x)$ is the interior normal to the boundary of $\Omega$.
Let $\psi $ be a smooth function in $\Omega$, such that $0\le\psi\le2$ in $\Omega$, $\psi=2$ in $\Omega\setminus \Omega^{\prime\prime}$,
$$
\psi(x) = d(x)\;\mbox{ in }\; \opr, \qquad  \frac{\partial\psi}{\partial\nu}\ge0 \;\mbox{ in }\; \Omega, \qquad \|\psi\|_{C^2(\Omega)}\le C(\Omega).
$$
Note we choose these just to ease some technicalities later, $\psi$ is equivalent to $d$ but is smooth, and we will consider $u/\psi$. \medskip

As in the proof of the interior weak Harnack inequality we set $k=\|f\|_{L^q(\Omega)}$ if $f\not\equiv0$, and let $k>0$ be arbitrary if $f=0$. Replace $u$ by $\tilde u=u+k$, which solves the same equation. We are going to show that, given $\varepsilon_0\in (0,1)$,
\begin{equation}\label{propineq}
\left( \int_{\Omega} \left(\frac{\tilde u}{\psi}\right)^s\right)^{1/s}\le C \left( \int_{\Omega} \left(\frac{\tilde u}{\psi}\right)^{\varepsilon_0}\right)^{1/\varepsilon_0}
\end{equation}
for each $s<1$. Here $C$  depends also on  $\varepsilon_0$. Proposition \ref{propgen} follows from \eqref{propineq}, \cite[Theorem 1.2]{Sir1}, and the standard inequalities: for $u,k, \alpha>0$,
$$
\min \{1,2^{\alpha-1}\}(u^\alpha + k^\alpha) \le (u+k)^\alpha\le \max \{1,2^{\alpha-1}\} (u^\alpha + k^\alpha),
$$
noting also that $\psi^{-s}$ is integrable in $\Omega$ for $s<1$, $\|\psi^{-s}\|_{L^1(\Omega)}=C(s,\Omega)$.\bigskip

We turn to the proof of \eqref{propineq}. Assume first $n>2$. Fix $s_0<1$ and let us prove \eqref{propineq} for $s=s_0$. Also fix some numbers $\sigma\in[0,1)$, $\gamma\in [\varepsilon_0, (1+s_0)/2]$, $r\in(\sigma\gamma,1)$. Later we will specify (and vary) these constants.

In all that follows $C$ will denote a constant which may vary from line to line and depends on the usual quantities, as well as on positive lower bounds on $1-s$, $1-\sigma$, $1-r$, $r-\sigma\gamma$.

We set,  for $m\in\mathbb{N}$,
$$\eta_m(x) = \left\{\begin{array}{ccc} md(x)&\mbox{if}& d(x)\le 1/m\\
1&\mbox{if}& d(x)\ge 1/m. \end{array} \right.
$$
We will often omit the subscript $m$, and write $\eta=\eta_m$.
Further denote
$$
v=v_\sigma = \psi^{-\sigma}\tilde u, \qquad\mbox{and}\qquad w=v^{\gamma/2}.
$$

We recall that $0<k\le \tilde u\le M<\infty$ in $\Omega$, with $M= \|u\|_{L^\infty(\Omega)} +k$.

It is easy to check that the function $\eta w\in H^1_0(\Omega)$ (near the boundary $\partial \Omega$ we have $\eta w\sim  d^{1-\frac{\sigma\gamma}{2}}$, and $\sigma\gamma<1$, so $\nabla(\eta w)$ is square-integrable up to the boundary). Our  goal is an estimate of the type
$$
\io \psi^{1+\gamma} |\nabla(\eta w)|^2 \le C\io (\eta w)^2 + \mathrm{negligible},
$$
 which together with Theorem \ref{wsob} will lead to a reverse H\"older inequality.\bigskip

We will use the test function
$$
\varphi= \psi^{1+r-\sigma\gamma} {\tilde u}^{\gamma-1} \eta^2 = \psi^{1-\sigma+r} v^{\gamma-1} \eta^2
$$
in the weak formulation of \eqref{first}
\begin{equation}\label{firstweak}
\io (A\nabla  u, \nabla\varphi) + \io b|\nabla  u|\varphi \ge -\io  f\varphi,
\end{equation}
valid for each $\varphi\in H^1_0(\Omega)$, $\varphi\ge0$ (by density), and will rewrite the resulting inequality in terms of $w$.

We compute, setting $ \tilde f=f/k$, $\|\tilde f\|_{L^q}\le1$, that for each $\epsilon>0$ we can find $C_\epsilon>0$ for which
\begin{eqnarray}
\left|\io f\varphi\right|&=&\left|\io \psi^{1-\sigma+r} \frac{f}{v} v^\gamma \eta^2\right| = \left|\io \psi^{1+r} \frac{f}{\tilde u} v^\gamma \eta^2\right|\nonumber\\
&\le& \io \psi^{1+r} \tilde f (\eta w)^2 \le \|\tilde f\|_{L^{q/2}}\,\|\psi^{1+r} (\eta w)^2\|_{L^{(q/2)^\prime}}\nonumber\\
&\le& \epsilon \|\psi^{1+r} (\eta w)^2\|_{L^{n/(n-2)}} +C_\epsilon \|\psi^{1+r} (\eta w)^2\|_{L^1}\nonumber\\
&\le& C\epsilon  \io \psi^{1+r}|\nabla(\eta w)|^2 + C_\epsilon\io (\eta w)^2,\label{compf}
\end{eqnarray}
In the third inequality we used that $q>n$, so $(q/2)^\prime<(n/2)^\prime=n/(n-2)$, H\"older and Young inequalities; to get the last inequality in \eqref{compf}, we applied Theorem \ref{wsob} with
$$
a=b=\frac{1+r}{2},\quad t = \frac{2n}{n-2}.
$$

\noindent {\it Remark.} For this computation we only need $f\in L^{q/2}$,  $q>n$.\bigskip

Further, since
\begin{equation}\label{grads}
\nabla u = \nabla \tilde u=\sigma \psi^{\sigma -1} v\nabla \psi  +\psi^\sigma \nabla v, \qquad v^{\gamma-1}\nabla v = \frac{1}{\gamma} \nabla w^2 = \frac{2}{\gamma}w \nabla w,
\end{equation}
we have, setting $\tilde b= |b.\nabla\psi|\in L^q(\Omega)$,
\begin{eqnarray*}
\left|\io b|\nabla u|\varphi\right| &\le & \io \tilde b \psi^r v^\gamma \eta^2 +
\frac{C}{\gamma}\io \psi^{1+r} \eta^2 w b|\nabla w|\\
&\le & \io \tilde b \psi^r  (\eta w)^2 + \frac{C}{\varepsilon_0} \io \psi^{1+r} \eta w b\left(|\nabla(\eta w)| + w|\nabla \eta|\right)\\
& =: &J_1 + (C/\varepsilon_0)|J_2| +(C/\varepsilon_0)|J_3|.
\end{eqnarray*}
We evaluate, by $q^\prime<n^\prime=n/(n-1)$, for every $\epsilon>0$,
\begin{eqnarray*}
J_1&\le& \|\tilde b\|_{L^q}\|  (\psi^{r/2}\eta w)^2\|_{L^{q^\prime}}\\
&\le & C\epsilon \|  \psi^{r/2}\eta w\|^2_{L^{2n/(n-1)}} +C_\epsilon  \|  (\psi^{r/2}\eta w)^2\|_{L^1}\\
&\le & C\epsilon \io \psi^{1+r}|\nabla(\eta w)|^2 + C_\epsilon\io  (\eta w)^2
\end{eqnarray*}
where in the last inequality we used Theorem \ref{wsob} with
$$a= \frac{1+r}{2}, \quad b=\frac{r}{2},\quad t= \frac{2n}{n-1}.
$$

We observe that
$$
\psi|\nabla \eta|\le \eta,
$$
so
$$
|J_3|\le C\io |b|\psi^r  (\eta w)^2,
$$
and $J_3$ can be evaluated exactly like $J_1$. On the other hand,
$$
|J_2|\le \epsilon \io \psi^{1+r}|\nabla(\eta w)|^2 +C_\epsilon \io \psi^{1+r} b^2(\eta w)^2,
$$
and the last integral can be evaluated exactly like  in \eqref{compf}, replacing $\tilde{f}$ there by $b^2\in L^{q/2}$.\medskip

In the following we denote with $\iok$ any integral such that for any $\epsilon>0$ there exists $C_\epsilon>0$ for which
 $$
| \iok|\le  \epsilon \io \psi^{1+r}|\nabla(\eta w)|^2 +C_\epsilon \io (\eta w)^2.
 $$
 We have just shown the second and the third term in \eqref{firstweak} have this property.\bigskip

We now turn to the highest order integral in \eqref{firstweak} - of $(A\nabla u,\nabla \varphi)$, where most care will be needed. We have
$$
\nabla\varphi = (\gamma-1)\psi^{1-\sigma+r}\eta^2v^{\gamma-2}\nabla v + (1-\sigma+r) \psi^{-\sigma +r} \eta^2v^{\gamma-1} \nabla \psi   + \psi^{1-\sigma+r}v^{\gamma-1} \nabla\eta^2,
$$
$$
v^{\gamma-2}(A\nabla v,\nabla v) = \frac{4}{\gamma^2}(A\nabla w, \nabla w),
$$
so, recalling \eqref{grads} and that $A$ is symmetric,
\begin{eqnarray}
\io (A\nabla u, \nabla\varphi)&=& -\frac{4(1-\gamma)}{\gamma^2} \io \psi^{1+r}\eta^2 (A\nabla w, \nabla w) + \frac{1}{\gamma}\io \psi^{1+r}(A\nabla w^2,\nabla \eta^2)\nonumber\\
&+& \frac{1+r+\sigma(\gamma-2)}{\gamma}\io \psi^r\eta^2(A\nabla w^2,\nabla \psi) + \sigma\io \psi^r w^2 (A\nabla \psi, \nabla \eta^2)\nonumber\\
&+& \sigma(1-\sigma + r) \io \psi^{r-1}w^2\eta^2 (A\nabla\psi , \nabla\psi)\nonumber\\
&=:&I_1+I_2+I_3+I_4+I_5.\label{thisone}
\end{eqnarray}

We have
$$
\io \psi^r \eta^2 (A\nabla w^2, \nabla \psi) = \io \psi^r (A\nabla (\eta w)^2, \nabla \psi) - \io \psi^r w^2 (A\nabla \psi,\nabla \eta^2),
$$
and, by the divergence theorem
\begin{eqnarray*}
 \io \psi^r (A\nabla (\eta w)^2, \nabla \psi) &=& -r\io \psi^{r-1}(\eta w)^2 (A\nabla \psi, \nabla \psi) + \io (A\nabla(\psi^r\eta^2w^2), \nabla \psi) \\
 &=&-r\io \psi^{r-1}(\eta w)^2 (A\nabla \psi, \nabla \psi) - \io \mathrm{div}(A\nabla\psi)\psi^r(\eta w)^2\\& +& \int_{\partial \Omega} \psi^r (\eta w)^2 (A\nabla\psi,-\nu).
 \end{eqnarray*}
 The integral on the boundary vanishes (recall the definition of $\psi$, as well as $w^2\le C\psi^{-\sigma\gamma}u^\gamma\le CM^\gamma\psi^{-\sigma\gamma}$ and $\sigma\gamma<r<1$), while the penultimate integral  can be evaluated exactly like $J_1$, since $\|\mathrm{div}(A\nabla\psi)\|_{L^q}\le C\Lambda$.

 We deduce that
\begin{eqnarray*}
I_3+I_4+I_5 &=& \left(\sigma(1-\sigma)-\frac{r}{\gamma}(1+r-2\sigma)\right)\io \psi^{r-1}(\eta w)^2 (A\nabla \psi, \nabla \psi) \\ &-&\frac{1+r-2\sigma}{\gamma}\io \psi^r w^2 (A\nabla \psi,\nabla \eta^2)+\iok.
\end{eqnarray*}

Next, to evaluate $I_2$ we observe that
$$
\io \psi^{1+r}(A\nabla w^2,\nabla \eta^2) = -(1+r) \io \psi^r w^2 (A\nabla \psi, \nabla \eta^2) + \io (A\nabla(\psi^{1+r}w^2),\nabla\eta^2),
$$
while by the divergence theorem
\begin{eqnarray*}
\io (A\nabla(\psi^{1+r}w^2),\nabla\eta^2)&=& \io (A\nabla\eta^2,\nabla(\psi^{1+r}w^2))\\
&=& \int_{\partial \Omega} \psi^{1+r} w^2 (A\nabla\eta^2,-\nu) - \io \mathrm{div}(A\nabla\eta^2)\psi^{1+r}w^2  \\
&=&  -\io A^\prime \nabla\eta^2\psi^{1+r}w^2 - \io \mathrm{tr}(AD^2\eta^2)\psi^{1+r}w^2,
\end{eqnarray*}
where $A^\prime$ is a matrix containing derivatives of the entries of $A$, so $A^\prime$ is bounded in $L^q$. Since $|\nabla\eta^2|\le C\psi^{-1}\eta^2$,  the integral on $\partial\Omega$ vanishes, and the first integral in the right-hand side of the last equality can again be evaluated like $J_1$. We will next deal with the last integral.

We fix a smooth orthonormal basis $(\tau,\nu)=(\tau_1,\ldots,\tau_{n-1},\nu)$, where $\nu(x)$ is the vector field we defined above, and let $T(x)$ be an orthogonal change-of-basis matrix between $(\tau(x),\nu(x))$ and $x$ (the $C^2$-norm of $T$ is bounded in terms of $\Omega$), so that $\nabla = \nabla_x=T\nabla_{\tau,\nu}$.
By the definition of $\eta=\eta_m$ we see that the directional derivative
$\partial_{\alpha\beta}(\eta^2)= 0$ for each pair of vectors $\alpha,\beta\in \{\tau_1,\ldots,\tau_{n-1},\nu\}$ such that $(\alpha,\beta)\not=(\nu,\nu)$, and
$$
\partial_{\nu\nu}(\eta^2)=2m^2\,\chi({\Omega_m^\prime}) - 2m\,\delta(\partial\Omega_m),
$$
where $\Omega_m^\prime=\{x\in\Omega\::\: \mathrm{dist}(x,\partial \Omega)<1/m\}$, $\Omega_m= \Omega\setminus\Omega_m^\prime$, $\chi$ denotes the characteristic function, and $\delta$ is the Dirac mass concentrated at $\partial\Omega_m$.

Therefore
\begin{eqnarray*}
\io (A\nabla(\psi^{1+r}w^2),\nabla\eta^2)&=& \left( -2m^2\int_{\Omega_m^\prime} + 2m \int_{\partial\Omega_m}\right)\tilde{a}_{nn}\psi^{1+r}w^2 + \iok
\end{eqnarray*}
where $\tilde{a}_{nn}$ denotes the last entry of the matrix $\tilde{A}=T^{-1}AT$, which has the property
\begin{equation}\label{tukuk}
(A\nabla\cdot, \nabla\cdot)=(\tilde{A}\nabla_{\tau,\nu}\cdot,\nabla_{\tau,\nu}\cdot).
\end{equation}

We compute
$$
0\le\lim_{m\to\infty} 2m \int_{\partial\Omega_m}\tilde{a}_{nn}\psi^{1+r}w^2 \le \lim_{m\to\infty} 2Cm\frac{1}{m^{1+r}}m^{\sigma\gamma}(\max_{\bar\Omega} \tilde{a}_{nn}) =0,
$$
since $\psi=1/m$ on $\partial\Omega_m$, $w^2\le CM^\gamma\psi^{-\sigma\gamma}$ (as above), and $r>\sigma\gamma$.

Finally,
\begin{eqnarray*}
I_2+\ldots+I_5&=&\frac{-2(1+r) + 2\sigma}{\gamma} \io \psi^r w^2 (A\nabla \psi, \nabla \eta^2)- \frac{2m^2}{\gamma}\int_{\Omega_m^\prime}\tilde{a}_{nn}\psi^{1+r}w^2\\
&+&\left(\sigma(1-\sigma)-\frac{r}{\gamma}(1+r-2\sigma)\right)\io \psi^{r-1}(\eta w)^2 (A\nabla \psi, \nabla \psi)+\iok + o(1),
\end{eqnarray*}
where $o(1)$ is a quantity which goes to zero as $m\to\infty$.

Therefore, by \eqref{thisone} and \eqref{firstweak},
\begin{eqnarray}
\frac{4(1-\gamma)}{\gamma^2} \io \psi^{1+r}\eta^2 (A\nabla w, \nabla w)&\le&\frac{-2(1+r) + 2\sigma}{\gamma} \io \psi^r w^2 (A\nabla \psi, \nabla \eta^2)\label{imeq1}\\ &+& \left(\sigma(1-\sigma)-\frac{r}{\gamma}(1+r-2\sigma)\right)\io \psi^{r-1}(\eta w)^2 (A\nabla \psi, \nabla \psi)\nonumber \\ &-& \frac{2m^2}{\gamma}\int_{\Omega_m^\prime}\tilde{a}_{nn}\psi^{1+r}w^2+\iok + o(1),\nonumber
\end{eqnarray}
For a moment we set in \eqref{imeq1}
$$
\sigma=0, \quad\mbox{i.e. }\; v=u,\; w=w_0=u^{\gamma/2}.
$$
Then all three constants in front of the integrals in the right-hand side of \eqref{imeq1} are negative, and $(1-\gamma)/\gamma^2$ is between two positive constants (by the assumption we made on $\gamma$), so by the uniform positivity of~$A$
\begin{equation}\label{tuk}
\io \psi^{1+r} \eta^2 |\nabla w_0|^2 \le \iok +o(1), \qquad\mbox{for all }\; r\in (0,1).
\end{equation}
We recall that the constant in $\iok$ in \eqref{tuk} is bounded by what we need, as long as $r$ is bounded away from $0$ and $1$.\bigskip

We go back to the general case $\sigma>0$. We start with the following trivial observation. By using the Young inequality on the terms involving $\tilde{a}_{\tau_i\nu}$ in the quadratic form $(\tilde{A}\xi, \xi)$, for any $\delta>0$ we can find $C_\delta>0$ such that for all $\xi\in\mathbb{R}^n$, $\xi=(\xi^\prime,\xi_n)$,
$$
\tilde{a}_{in}\xi_{i}\xi_n\ge -(\delta/(\lambda n^2)) \xi_n^2  - C_\delta \xi_{i}^2\ge -(\delta/n^2) \tilde{a}_{nn}\xi_n^2 - C_\delta |\xi^\prime|^2,
$$
and writing
$$(\tilde{A}\xi,\xi) = (\tilde{A}^\prime\xi^\prime,\xi^\prime) + \tilde{a}_{nn} \xi_n^2 + \sum_{i=1}^{n-1} \tilde{a}_{in}\xi_{i}\xi_n,
$$
we have
\begin{equation}\label{imeq2}
(1-\delta)\tilde{a}_{nn} \xi_n^2\le (\tilde{A}\xi,\xi) +C_\delta |\xi^\prime|^2.
\end{equation}

From now on we set
$$
r=\gamma.
$$

By \eqref{tukuk}, plugging \eqref{imeq2} into \eqref{imeq1},
\begin{eqnarray}
(1-\delta) \io \psi^{1+\gamma}\eta^2 \tilde{a}_{nn}(\partial_\nu w)^2 &\le& C_\delta \io \psi^{1+\gamma}\eta^2 |\nabla_\tau w|^2 -  D_1 \io \psi^\gamma w^2 (A\nabla \psi, \nabla \eta^2)\nonumber\\ &+& D_2\io \psi^{\gamma-1}(\eta w)^2 (A\nabla \psi, \nabla \psi)\label{lpoc}\\ &-& D_3m^2\int_{\Omega_m^\prime}\tilde{a}_{nn}\psi^{1+\gamma}w^2+\iok + o(1),\nonumber
\end{eqnarray}
where
$$
D_1=\frac{\gamma(1+\gamma-\sigma)}{2(1-\gamma)}\quad
D_2=\frac{\gamma^2\left(-\sigma^2 + 3\sigma -1-\gamma\right)}{4(1-\gamma)},\quad
D_3=\frac{\gamma}{2(1-\gamma)}.
$$

However, by using $w^2=w_\sigma^2=u^\gamma/\psi^{\sigma\gamma}$ and $\nabla_\tau \psi=0$ in $\opr$, $\psi\ge1$ in $\Omega\setminus\opr$,
\begin{equation}\label{checks}
\io \psi^{1+\gamma}\eta^2 |\nabla_\tau w|^2 \le C \io \psi^{1+\gamma-\sigma\gamma} \eta^2 (w_0^2+|\nabla w_0|^2) \le \iok +o(1);
\end{equation}
the second inequality follows from \eqref{compf} and \eqref{tuk} with $r$ replaced by $(1-\sigma)\gamma$.

We have $\psi^{-1}\le 1$ in $\Omega\setminus\opr$ and $(A\nabla \psi, \nabla \psi)=\tilde{a}_{nn}$ in $\opr$,
so
\begin{eqnarray}
\io \psi^{\gamma-1}(\eta w)^2 (A\nabla \psi, \nabla \psi)&\le& \io \psi^{\gamma}(\eta w)^2 (A\nabla \psi, \nabla \psi) + \int_{\opr} \psi^{\gamma-1}(\eta w)^2 \tilde{a}_{nn}\nonumber\\ &\le&  \iok + \int_{\opr}  \psi^{\gamma-1}(\eta w)^2 \tilde{a}_{nn},\label{firin}
\end{eqnarray}
where $\iok$ is evaluated exactly like $J_1$ above. Note $|(A\nabla \psi, \nabla \psi)|\le C$.

We further compute, writing $\partial=\partial_\nu$,
\begin{eqnarray*}
\gamma\io \psi^{\gamma-1} w^2 \eta^2 \tilde{a}_{nn} \partial\psi  &+& \io \psi^\gamma\partial (w^2) \eta^2 \tilde{a}_{nn}  \\
&=& \io \partial(\psi^\gamma w^2) \eta^2 \tilde{a}_{nn} \\
=-\io \psi^\gamma w^2\partial(\eta^2) \tilde{a}_{nn} &-& \io \psi^\gamma w^2\eta^2 \partial \tilde{a}_{nn}  + \int_{\partial\Omega} + \iok\\
=-\io \psi^\gamma w^2\partial(\eta^2) \tilde{a}_{nn} &+& \iok
\end{eqnarray*}
(the boundary term again vanishes).

Hence (recall $\partial \psi\ge0$ in $\Omega$ and $\partial \psi = 1$ in $\opr$)
\begin{eqnarray}
\int_{\opr}  \psi^{\gamma-1}(\eta w)^2 \tilde{a}_{nn}&\le & \io \psi^{\gamma-1} w^2 \eta^2 \tilde{a}_{nn} \partial\psi\label{secin}\\
&=& -\frac{1}{\gamma} \io \psi^\gamma\partial (w^2) \eta^2 \tilde{a}_{nn} -\frac{1}{\gamma}\io \psi^\gamma w^2\partial(\eta^2) \tilde{a}_{nn} +\iok\nonumber
\end{eqnarray}

Further, we have
\begin{eqnarray}
\left|\io \psi^\gamma \partial (w^2) \eta^2 \tilde{a}_{nn}\right| &\le& 2\io (\psi^{(\gamma+1)/2}|\partial w| \eta \tilde{a}_{nn}^{1/2})(\psi^{(\gamma-1)/2} w \eta \tilde{a}_{nn}^{1/2})\nonumber\\
&\le& \frac{2}{\gamma} \io \psi^{\gamma+1}(\partial w)^2 \eta^2 \tilde{a}_{nn}  + \frac{\gamma}{2} \io \psi^{\gamma-1}w^2 \eta^2\tilde{a}_{nn}\label{thiin}\\
&=& \frac{2}{\gamma} \io \psi^{\gamma+1}(\partial w)^2 \eta^2 \tilde{a}_{nn} + \frac{\gamma}{2} \int_{\opr} \psi^{\gamma-1}w^2 \eta^2\tilde{a}_{nn} + \iok\nonumber
\end{eqnarray}
since again $\psi^{-1}\le 1$ in $\Omega\setminus\opr$.

Combining \eqref{secin} with \eqref{thiin}, we get from \eqref{firin}
$$
\io \psi^{\gamma-1}(\eta w)^2 (A\nabla \psi, \nabla \psi) \le \frac{4}{\gamma^2} \io \psi^{\gamma+1}(\partial w)^2 \eta^2 \tilde{a}_{nn} - \frac{2}{\gamma}\io \psi^\gamma w^2\partial(\eta^2) \tilde{a}_{nn} +\iok.
$$
Hence, noticing that $(A\nabla \psi, \nabla \eta^2)=\partial(\eta^2) \tilde{a}_{nn}$, we infer from \eqref{lpoc} and \eqref{checks}
\begin{eqnarray}
(1-\delta) \io \psi^{1+\gamma}\eta^2 \tilde{a}_{nn}(\partial w)^2 &\le& -(D_1+\frac{2D_2}{\gamma})  \io \psi^\gamma w^2\partial(\eta^2) \tilde{a}_{nn}\nonumber\\ &+& \frac{4D_2}{\gamma^2} \io \psi^{\gamma+1}(\partial w)^2 \eta^2 \tilde{a}_{nn}\label{thatone}\\ &-& D_3m^2\int_{\Omega_m^\prime}\tilde{a}_{nn}\psi^{1+\gamma}w^2+\iok + o(1),\nonumber
\end{eqnarray}
Now $\sigma\in (0,1)$ implies
\begin{equation}\label{fine1}
\sigma^2-3\sigma+2>0
\end{equation}
 which in turn guarantees precisely  that
\begin{equation}\label{fine2}
\frac{4D_2}{\gamma^2}<1.
\end{equation}
Hence, by choosing $\delta>0$ sufficiently small, the second term in the right-hand side of \eqref{thatone} can be absorbed in the left-hand side.

It is easy to check that
\begin{equation}\label{fine3}
D_1+\frac{2D_2}{\gamma}=\gamma\frac{(2-\sigma)\sigma}{2(1-\gamma)}\ge0.
\end{equation}

Recalling \eqref{checks} and \eqref{thatone}, we have thus shown that for some uniformly positive constants $\lambda_0, d_0$,
$$
\lambda_0\io \psi^{1+\gamma} \eta^2 |\nabla w|^2\le -d_0m^2\int_{\Omega_m^\prime}\psi^{1+\gamma}w^2+\iok + o(1).
$$
We can assume $\lambda_0<d_0$, by further diminishing $\lambda_0$ if necessary.

We have $|\nabla\eta|^2=m^2\chi(\Omega_m^\prime)$. Hence
\begin{eqnarray*}
\frac{\lambda_0}{2}\io \psi^{1+\gamma}  |\nabla(\eta\psi)|^2&\le& \lambda_0 \io \psi^{1+\gamma} |\nabla\eta|^2 w^2+\lambda_0\io \psi^{1+\gamma} \eta^2 |\nabla w|^2+\iok + o(1)\\
&\le&(\lambda_0-d_0)m^2\int_{\Omega_m^\prime}\psi^{1+\gamma}w^2+\iok + o(1)\\
&\le&\iok + o(1).
\end{eqnarray*}

Taking $\epsilon = \lambda_0/4$ in the definition of $\iok$, we get
$$
\io \psi^{1+\gamma} |\nabla(\eta w)|^2 \le C\io (\eta w)^2 + o(1),
$$
so, by Theorem \ref{wsob} applied with $a=(1+\gamma)/2$, $b=0$, we obtain
$$
\|(\eta w)^2\|_{L^{\rho}(\Omega)}\le C\|(\eta w)^2\|_{L^{1}(\Omega)} +o(1),
$$
where
$$
\rho=\frac{t}{2}=\frac{n}{n-1+\gamma}.
$$

Letting $m\to\infty$, by the definition of $w$, $\eta=\eta_m\nearrow1$ and the monotone convergence theorem, we get
$$
\left\|\frac{\tilde u}{\psi^\sigma}\right\|_{L^{\rho\gamma}(\Omega)}\le C\left\|\frac{\tilde u}{\psi^\sigma}\right\|_{L^{\gamma}(\Omega)},
$$
as long as the right-hand side is finite. This is so in particular if $\gamma= \varepsilon_0$ (recall $\sigma\le1$), since we already know that $\left(\frac{\tilde u}{\psi}\right)^{\varepsilon_0}\in L^1$, by \cite[Theorem 1.2]{Sir1}. Set $a_1:=\varepsilon_0$.

Hence $\frac{u}{\psi^\sigma}\in L^{a_2}$, $a_2=\frac{n}{n-1+a_1}a_1$. Taking in the above argument $\gamma = \frac{n}{n-1+a_k}a_k$ results in the iteration
$$
\frac{\tilde u}{\psi^\sigma}\in L^{a_k},\qquad \left\|\frac{\tilde u}{\psi^\sigma}\right\|_{L^{a_{k}}(\Omega)}\le C^{k-1} \left\|\frac{\tilde u}{\psi^\sigma}\right\|_{L^{a_{1}}(\Omega)},
$$
\begin{equation}\label{recur}
a_1=\varepsilon_0, \qquad a_{k+1} =\frac{n}{n-1+a_k}a_k, \quad k\in\mathbb{N},
\end{equation}
It is easy to check that the recursively defined  sequence $\{a_k\}$ is increasing as long as $a_1=\varepsilon_0<1$, and
$$
 \lim_{k\to\infty} a_k = 1.
$$
We claim that the proof is finished after a finite number $k_0$ of iterations, where $k_0$ is the first index such that
$$
a_{k_0}\ge \frac{1+s_0}{2} = \frac{s_0}{\sigma},
$$
where the latter equality is how we make our overall choice of $\sigma$,
$$
\sigma:=\frac{2s_0}{1+s_0}.
$$
To prove the claim, by the H\"older inequality and $s={s_0}<1$,
\begin{eqnarray*}
\io \left(\frac{\tilde u}{\psi}\right)^{s} &=& \io \frac{\tilde u^{s}}{\psi^{\sigma {s}}}
\frac{1}{\psi^{(1-\sigma)s}}\\
&\le& \left(\io \frac{\tilde u^{s/\sigma}}{\psi^s}\right)^\sigma \left( \io \frac{1}{\psi^s}\right)^{1-\sigma}\\
&\le& C\left\|\frac{\tilde u}{\psi^\sigma}\right\|_{L^{s/\sigma}(\Omega)}^s \le C\left\|\frac{\tilde u}{\psi^\sigma}\right\|_{L^{a_{k_0}}(\Omega)}^s \\
&\le&   C^{k_0}\left\|\frac{\tilde u}{\psi^\sigma}\right\|_{L^{a_1}(\Omega)}^s \le C\left\|\frac{\tilde u}{\psi}\right\|_{L^{\varepsilon_0}(\Omega)}^s,
\end{eqnarray*}
and \eqref{propineq} is proved. \hfill $\Box$\bigskip

Finally, if $n=2$, we use the last part of Theorem \ref{wsob}. Repeating the above, with trivial modifications, we can show that
$$
\io \psi^{1+\gamma} |\nabla(\eta w)|^\alpha \le C_\alpha\io (\eta w)^2 + o(1),
$$
for each $\alpha<2$, with $C_\alpha$ bounded in terms of a lower bound on $2-\alpha$. Exactly as above we can set up an iteration process which produces a sequence $a_k$ converging to a number  $a(\alpha)<1$. We readily see that $a(\alpha)\to1$ as $\alpha\to2$. So for each initially fixed $s_0<1$ we can choose $\alpha<2$ such that $a(\alpha)>s_0$, and a finite number of iterations give a bound for $u/\psi$ in $L^{s_0}$. The technical details are left to the interested reader.

Proposition \ref{propgen} and Theorem \ref{bwhistrong} are proved.

\subsection{Proofs of Theorem \ref{bwhiweak} and Theorem \ref{bwhigrad}}

By the same argument as in the beginning of the previous section, Theorem~\ref{bwhiweak} and Theorem \ref{bwhigrad} have local versions around any point on $\partial\Omega$, so it is sufficient to prove them under the hypothesis of Proposition \ref{propgen}.

Theorem \ref{bwhiweak} is  simpler than Theorem \ref{bwhistrong}, and follows only from \eqref{tuk}. This inequality together with Theorem \ref{wsob} applied with $a=(1+r)/2$, $b=0$, implies that for each $\gamma<1$, $r<1$,
\begin{equation}\label{iterr}
\left|\tilde u\right|_{L^{\rho\gamma}(\Omega)}\le C\left|\tilde u\right|_{L^{\gamma}(\Omega)},\qquad \rho=\frac{n}{n-1+r}
\end{equation}
as long as the right-hand side is finite.

We can assume $p\ge1$. Fix $r_0<1$ so small that
$$
\frac{n}{n-1+r_0}= \frac{1}{2}\left( p+ \frac{n}{n-1}\right)
$$
Fix $k_0$ such that
$$
\left(\frac{n}{n-1+r_0}\right)^{k_0}\varepsilon_0>1
$$
and then $\gamma_0<\varepsilon_0$ such that
$$
\left(\frac{n}{n-1+r_0}\right)^{k_0}\gamma_0=\delta_0,
$$
where $\delta_0<1$ is so close to $1$ that
$$
\frac{n}{n-1+r_0}\,\delta_0=\frac{\delta_0}{2}\left( p+ \frac{n}{n-1}\right)>p.
$$
With these choices, $k_0+1$ iterations of \eqref{iterr} starting from $\gamma_0$ give Theorem~\ref{bwhiweak}.
\bigskip

We now prove Theorem \ref{bwhigrad}. Fix $s\in(0,1)$. We use the following well-known inequality, valid for bounded weak Sobolev supersolutions of \eqref{first}. For each $\eta\in  C_{\mathrm{c}}^1(\Omega)$, by testing \eqref{first} with $\eta^2 {\tilde u}^{s-1}$,
\begin{equation}\label{gilg}
\int_\Omega \eta^2{\tilde u}^{s-2}|\nabla u|^2\le C\int_\Omega (\eta^2 + |\nabla \eta|^2){\tilde u}^{s},
\end{equation}
where $C$ is bounded in terms of positive lower and upper bounds for $1-s$. This follows from the computation on pages 195-196 in \cite{GT}, in particular inequalities (8.52)-(8.53)-(8.54) with $\beta=s-1$ there.

By density, the same inequality is valid for any $\eta\in H^1_0(\Omega)$. Since $d^{r}$ has square-integrable gradient  for each $r>1/2$, we can apply \eqref{gilg} with $$\eta=\psi^{1-\frac{s}{2}}\le Cd^{1-\frac{s}{2}}.$$

Then
$$
\eta^2 + |\nabla\eta|^2\le Cd^{-s},
$$
so \eqref{gilg} becomes
$$
\int_\Omega \eta^2{\tilde u}^{s-2}|\nabla u|^2\le C\int_\Omega \left( \frac{\tilde u}{d}\right)^s,
$$
and we obtain by applying the H\"older inequality
\begin{eqnarray*}
\int_\Omega |\nabla u|^s &=& \int_\Omega \left(\eta^{{s}}\, \tilde u^{-\frac{s(2-s)}{2}}\,|\nabla u|^s\right) \left(\eta^{-{s}}\, \tilde u^{\frac{s(2-s)}{2}}\right)\\
&\le& \left(\int_\Omega \eta^2{\tilde u}^{s-2}|\nabla u|^2\right)^{s/2}\:\left(\int_\Omega \eta^{-\frac{2s}{2-s}} \tilde u^s\right)^{(2-s)/2}\\
&\le& C\int_\Omega \left( \frac{\tilde u}{d}\right)^s\le C\left(\int_\Omega \left( \frac{ u}{d}\right)^s + \|f\|_{L^q(\Omega)}^s\right).
\end{eqnarray*}
We conclude the proof of Theorem \ref{bwhigrad} with the help of Theorem \ref{bwhistrong}.\hfill $\Box$ \bigskip

In the end we recall, for completeness, that in case $u$ is a solution, rather than just a supersolution, the gradient of $u$ is bounded pointwise by the quantity $u/d$.  Indeed, by standard elliptic regularity we know that $u\in C^{1,\alpha}_{\mathrm{loc}}(\Omega)$, with the gradient estimate
$$
\sup_K|\nabla u|\le C \left(\sup_{K^\prime} u + \|f\|_{L^q(K^\prime)}\right),
$$
for each $K\subset\subset K^\prime \subset\subset \Omega$, with $C$ depending of course on $K$, $K^\prime$. Hence by the Harnack inequality
\begin{equation}\label{gradharn}
\sup_K|\nabla u|\le C \left(\inf_{K^\prime} u + \|f\|_{L^q(\Omega)}\right).
\end{equation}
 Fix $x_0\in \Omega^\prime$ and $d=d(x_0)$. We apply \eqref{gradharn} to the function $\tilde u(x)=u(x_0+dx)$, which satisfies the same equation with $b$ replaced by $db$ and $f$ replaced by $d^2f$ (but $d\le 1$), and with $K=B_{1/2}(0)$, $K^\prime=B_{3/4}(0)$, $\Omega=B_1(0)$. We deduce
$$
 d(x_0)|\nabla u(x_0)|=|\nabla \tilde u(0)|\le C \left(\tilde u(0)+ \|f\|_{L^q(\Omega)}\right)= C \left(u(x_0)+ \|f\|_{L^q(\Omega)}\right).
 $$
By \eqref{gradharn} with $K=\Omega\setminus\Omega^\prime$, the same  is valid for any $x_0\in \Omega$ with $d(x_0)\ge1$. Thus
$$
|\nabla u|^s \le C\left(\left|\frac{u}{d}\right|^s + \frac{\|f\|_{L^q(\Omega)}^s}{d^s}\right)\qquad \mbox{in }\;\Omega,
$$which is another way to infer Theorem \ref{bwhigrad}  from Theorem \ref{bwhistrong}, if we have a solution.\bigskip

\noindent {\it Remark. } The last argument, combined with \cite[Theorem 1.2]{Sir1}, also shows that a nonnegative solution $u$ of a general non-divergence form inequality as in \cite{Sir1}  is such that $|\nabla u|^\varepsilon\in L^1(\Omega)$, for some $\varepsilon>0$. This observation complements the results in \cite{Sir1}.

\section{Appendix}
We consider the inequality \eqref{first} under the hypotheses \eqref{hypocoef}.
We first give a detailed statement and explanation why Theorem 1.2 from \cite{Sir1} applies to the type of supersolutions we consider here.

\begin{thm}\label{bwhieps} Let $A,b,f$ be as in \eqref{hypocoef} and $u\in H^1_{\mathrm{loc}}(\Omega)\cap L^\infty(\Omega)$  such that
\begin{equation}\label{defsupersolapp}
\int_\Omega A\nabla u.\nabla \varphi +\int_\Omega\varphi b|\nabla u| \ge -\int_\Omega f\varphi,\qquad \mbox{for\ each}\; \varphi\in C^\infty_0(\Omega), \varphi\ge0.
\end{equation}
Then there exist $\varepsilon>0$ depending on $n$, $\lambda$, $\Lambda$, $q$, $\|b\|_{L^q(\Omega)}$, and $C>0$ depending  on $n$, $\lambda$, $\Lambda$, $q$, $\|b\|_{L^q(\Omega)}$,   $\mathrm{diam}(\Omega)$, the $C^{1,1}$-norm of $\partial\Omega$, such that
 \begin{equation}\label{ineqapp}
\left( \int_{\Omega} \left(\frac{u}{d}\right)^\varepsilon\right)^{1/\varepsilon}\le C\left(\inf_{\Omega} \frac{u}{d} +\|f\|_{L^q(\Omega)}\right).
 \end{equation}
\end{thm}

\noindent{\it Proof.} For $H^1$-weak (super)solutions in the sense of \eqref{defsupersolapp} the following facts are known. They can be found in \cite[Chapter 8]{GT} and in \cite{Tr7}, \cite{LU} for operators with more general (unbounded) coefficients.
\begin{enumerate}
\item weak maximum/comparison principle -- see \cite[Theorem 8.1]{GT};
\item generalized maximum principle -- see \cite[Theorem 8.16]{GT};
\item solvability of the Dirichlet problem and global $C^{1,\alpha}$-estimates -- see \cite[Theorem 8.3]{GT}, \cite[Chapter 8.11]{GT} ;
\item interior Harnack inequality -- see  \cite[Theorem 8.18]{GT}.
\end{enumerate}

We can repeat practically verbatim the proof in \cite[Section 4]{Sir1}, replacing by 1.-4. above the results from the theory of viscosity solutions used in \cite{Sir1}. More precisely, that Theorem \ref{bwhieps} above follows from the growth lemma  \cite[Theorem 4.2]{Sir1} uses only measure theory and the interior Harnack inequality. To prove  \cite[Theorem 4.2]{Sir1} for supersolutions in the sense of \eqref{defsupersolapp} we repeat the proof on pages 7479-7480 in \cite{Sir1}. In that proof we replace the use of \cite[Theorem 2.1]{Sir1} by a rescaled version of \cite[Theorem 8.16]{GT} in a domain with width $\delta$, observing that \cite[Theorem 2.3.2]{Sir1} can be proved in exactly the same way in our setting here, by comparison, solvability of the Dirichlet problem, and $C^1$-estimates. We also observe that if a regular function $u\in W^{2,q}_{\mathrm{loc}}(\Omega)$ is a supersolution in the sense of \eqref{defsupersolapp} then it satisfies the elliptic  inequality \eqref{first} almost everywhere and after removing the divergence in the operator $u$ is also a strong solution of
$$
\mathcal{M}^-_{\lambda, \Lambda}(D^2u) -\tilde b |Du| \le \mathrm{tr}(AD^2u) -\tilde b |Du|\le  f
$$
in $\Omega$, for some $\tilde b$ such that $\|\tilde b\|_{L^q(\Omega)}\le \|b\|_{L^q(\Omega)}+ \sum\|Da_{ij}\|_{L^q(\Omega)}$.
\hfill $\Box$\bigskip

Next we  record some essentially known facts about various types of weak supersolutions, which in particular imply that all notions of supersolutions we recalled after Definition \ref{def1} are included in that definition. The following is also an alternative way to see that \cite[Theorem 1.2]{Sir1} applies to the type of supersolutions we consider here.

For precise definitions and the general theory of $C$-viscosity solutions we  refer to \cite{CIL}, for   $L^p$-viscosity solutions to \cite{CCKS}.
These two notions are coherent with respect to each other, specifically, if a continuous function is a $L^p$-viscosity solution of an equation with continuous coefficients then it is a $C$-viscosity solution (this is obvious by the definitions), and vice versa (by \cite[Proposition 2.9]{CCKS}).

\begin{prop}\label{prop1}
Viscosity solutions of \eqref{first} are $H^1$-weak supersolutions, i.e. \eqref{defsupersolapp} holds.
\end{prop}

\noindent{\it Proof}. This follows from the more general fact that any locally bounded function $w$ which satisfies the comparison principle with respect to regular subsolutions (i.e. if $v\in W^{2,q}$ is a subsolution and  $w\ge v$ on the boundary of a domain then $w\ge v$ in the domain -- note this is a basic property of viscosity supersolutions)  belongs to $H^1_{\mathrm{loc}}(\Omega)$ and is a weak Sobolev supersolution.  See for instance
\cite[Theorem 2]{HH}, and \cite[Proposition 14]{HH}. For more general $p$-laplacian like operators we refer to the book \cite{HKM}. \hfill $\Box$

\begin{prop}\label{prop2} If $v\in H^1_{\mathrm{loc}}(\Omega)\cap L^\infty(\Omega)$ satisfies  \eqref{defsupersolapp}, where $A(x)$ is a bounded uniformly positive matrix in $\Omega$ (for this we do not need any regularity for $A$) and  $b\in L^q_{\mathrm{loc}}(\Omega)$, $f\in L^{q/2}_{\mathrm{loc}}(\Omega)$, $q>n$, then $v$ is a lower semi-continuous function (after redefinition on a set of measure zero) which satisfies the definition of a viscosity supersolution.
\end{prop}

\noindent{\it Proof.}
Given $B_{2r_0}=B_{2r_0}(x_0)\subset \Omega$ and $r<r_0$ we define
$$
m(r)= \inf_{B_r} v
$$
(we recall that inf stands for essential infimum). By the interior weak Harnack inequality (Theorem 8.18 in GT) we have
\begin{eqnarray*}
0&\le& \frac{1}{|B_{2r}|}\int_{B_{2r}} (v(x)-m(2r)) \,dx\\
&\le& C(m(r)-m(2r)) + Cr^{2(1-n/q)}\|f\|_{L^{q/2}(B_{r_0})}.
\end{eqnarray*}
Since $m(r)$ is bounded and monotone, and $q>n$, the latter quantity tends to zero as $r\to 0$. Hence
$$
\liminf_{x\to x_0} v(x) = \lim_{r\to0} m(2r) = \lim_{r\to 0} \frac{1}{|B_{2r}|}\int_{B_{2r}(x_0)} v(x) \,dx.
$$
The last limit is $v(x_0)$ for almost every $x_0$, by the Lebesgue differentiation theorem. But the quantity $\liminf_{x\to x_0} v(x) $ is always lower semi-continuous in $x_0$, for any $v$. \bigskip

Assume now that $v$ does not satisfy the definition of a $C$-viscosity or $L^{q/2}$-viscosity supersolution of \eqref{first} in $\Omega$. This means  there exist a ball $B=B_{2r_0}(x_0)\subset\Omega $ and a function $\psi\in W^{2,q/2}(B)$ (note $W^{2,q/2}(B)\subset H^1(B)\cap C(\overline{B}))$ such that $\psi$ touches $v$ from below at $x_0$, but for some $\delta>0$
$$
-L\psi \le f-\delta \qquad \mbox{a. e.  in}\; B.
$$
Let $\theta\in H^1_0(B)$ be the unique solution (see \cite[Theorem 8.3]{GT}, and \cite{Tr7} where equations with unbounded coefficients are treated) of
$$
-L\theta = \delta\quad \mbox{in}\; B, \qquad \theta=0 \quad \mbox{on}\; \partial B.
$$
By De Giorgi's classical result (\cite[Theorem 8.24]{GT}) or by regularity we know that $\theta$ is continuous in $B$. By the weak maximum principle (\cite[Theorem 8.1]{GT}) and the interior Harnack inequality (\cite[Theorem 8.18]{GT}) we have $\theta>0$ in $B=B_{2r_0}$. Hence $\theta\ge \theta_0>0$ in $B_{r_0}$ for some positive constant $\theta_0$.

Thus, the function $w=v-\psi-\theta$ is such that $w(x_0)\le -\theta_0$ and $w$ satisfies in the weak Sobolev sense
$$
-Lw\ge0 \quad \mbox{in}\; B, \qquad w\ge 0 \quad \mbox{on}\; \partial B.
$$
By the maximum principle (\cite[Theorem 8.1]{GT}) $w\ge0$ in $B$, a contradiction.\bigskip

\noindent{\it Remark.} For possible further reference we note that in the definition of a viscosity solution the minimum at $x_0$ of $v-\psi$ can be assumed to be strict. For a full proof of this fact for equations with unbounded ingredients we refer to \cite[Lemma 2.10]{No}. \smallskip

\begin{prop} Potential theory supersolutions coincide with viscosity supersolutions.
\end{prop}

\noindent{\it Proof.} This is a simple exercise, using the function $\theta$ from the proof of Proposition \ref{prop2}. If a function $v$ is not a viscosity supersolution, then as in the previous proof $\psi + \theta$ is a strong solution of $-L(\psi + \theta)\le f$ in $B$, $\psi + \theta\le v$ on $\partial B$, so by the comparison principle $\psi +\theta\le z$, where $z$ is the strong solution of $-Lz= f$ in $B$, $z= v$ on $\partial B$. If $v$ were a potential theory supersolution we would have $z\le v$, leading to the same contradiction as in the proof of Proposition \ref{prop2}, since $\psi$ touches $v$ from below  in $B$.

Conversely, if $v$ is not a potential theory supersolution, there is a ball $B$ such that for $-Lz= f$ in $B$, $z= v$ on $\partial B$ but $z>v$ somewhere in $B$. Then for sufficiently small $\delta>0$ and some $C>0$ the function $z-\theta - C$ touches $v$ from below in $B$ and solves $-L(z-\theta - C)=f-\delta<f$, a contradiction with the definition of a viscosity supersolution. \hfill $\Box$

\end{document}